\documentclass[12pt]{article}

\usepackage{layout, graphicx, amsmath, amsthm, amssymb, euscript, enumerate, bbold, bbm, graphicx, fancyhdr,xcolor,stmaryrd,mathrsfs,pdfsync}

\usepackage{subcaption}
\usepackage{capt-of}% or \usepackage{caption} Para la tabla/fig

\makeatletter
\def\th@newremark{\th@remark\thm@headfont{\bfseries}}
\makeatletter

\theoremstyle{newremark}
\newtheorem{theorem}{Theorem}[section]

\newtheorem{proposition}[theorem]{Proposition}

\newtheorem{example}[theorem]{Example}

\newcommand{\N}{\mathbb{N}}

\newcommand{\Z}{\mathbb{Z}}

\begin{document}
	\title{Coalescent point process of branching trees in a varying environment}
	\author{  Airam~Blancas\footnote{Department of Statistics, ITAM, Mexico.
			Email: {airam.blancas@itam.mx}}
		\and 
		Sandra~Palau\footnote{Department of Probability and Statistics, IIMAS, UNAM, Mexico. Email: {sandra@sigma.iimas.unam.mx} }}
	\maketitle
	\vspace{-.2in}
	
	\begin{abstract}
Consider an arbitrary large population at the present time, originated at an unspecified arbitrary large time in the past, where individuals within the same generation independently reproduce forward in time, sharing a common offspring distribution that may vary across generations. In other words, the reproduction is  driven by a Galton-Watson process in a varying environment. The genealogy of the current generation, traced backward in time, is uniquely determined by the coalescent point process $(A_i, i\geq 1)$, where $A_i$ denotes the coalescent time between individuals $i$ and $i+1$. In general, this process lacks the Markov property. In constant environment, Lambert and Popovic (2013) proposed a Markov process of point measures to reconstruct the coalescent point process. We provide a counterexample showing that their process lacks the Markov property. The main contribution of this work is to propose a vector valued Markov process $(B_i,i\geq 1)$, that can reconstruct the genealogy, with finite information for every $i$. Additionally, in the case of linear fractional offspring distributions, we establish that the variables of the coalescent point process $(A_i, i\geq 1)$ are independent and identically distributed.\\

\textbf{MSC2020 subject classifications:} 60J80; 60J10; 60J85; 92D25.
	\end{abstract}

\section{Introduction}
Lambert and Popovic \cite{amaurylea2013} studied the backward genealogy of a random population where the forward dynamics are produced by a branching process, in either discrete or continuous-state space. In their model, the population at the present time could be arbitrarily large, originating from an unspecified distant past. In the discrete state case, they employed a monotone planar embedding tree where lines of descent do not intersect. The $i$-th individual in the past $n$-th generation, is represented by $(n,i)$ for any $n\in \Z_{-}:=\{0,-1,-2,\dots\}$ and $i\in \N:=\{1,2,\dots\}$. Individuals at the present generation are simply denoted  by $i$ instead of $(0,i)$. 

They demonstrated that the genealogy of the current generation, traced backward in time, is uniquely determined by the so-called \textit{coalescent point process}, denoted as $(A_i,i\geq 1)$, where $A_i$ represents the coalescent time between individuals $i$ and $i+1$. Generally, this process is not Markov, and characterizing its law poses challenges. To address this issue, they constructed a sequence-valued process $(D_i, i\geq 1)$ such that for every $i$, $A_i$ is the first non-zero entry of $D_i$. Their construction proceeds as follows: 
we follow its ancestral line (the so called spine). For every $n\in \N$, we consider the subtree attached to the spine with root at generation $-n$. We denote by $D_i(n)$ the number of children of the root at the right hand side of the spine, whose descendants remain alive at the present generation. Although this process is Markov, $D_i$ and $D_{i+1}$ consist of the same infinite sequence for every $i$, except for a finite number of entries. Consequently, $(D_i, i\geq 1)$ contains considerable repetitive information. Thus, they constructed a process $( \widetilde B_i, i\geq 1)$ by removing \textit{some} information from $(D_i, i\geq 1)$. They claimed that $(\widetilde B_i, i\geq 1)$ contains the minimal amount of information needed to construct $(A_i, i\geq 1)$ while remaining Markov. However, there was a mistake in their proof, and as exhibited in Example \ref{counterexample}, $(\widetilde B_i, i\geq 1)$ is not always Markov.

In this paper, our objective is to define a coalescent point process for a population driven by a Galton-Watson process in a varying environment. Here, individuals within the same generation share the same offspring distribution, but these distributions vary across generations. Our main contribution is to propose a new Markov process $(B_i, i\geq 0)$ capable of reconstructing the genealogy of a population in a varying environment, using finite information for every $i$. This property plays an important role, particularly in real-life applications, where evolutionary biologists seek to reconstruct the genealogy of a population based on a finite sample of individuals. For such purposes, employing the Markov process $(D_i, i\geq 1)$ as a model is impractical, as it necessitates an infinite sequence of values. Notice that in the case of a constant environment, $(B_i, i\geq 0)$ is a Markov process with a finite amount of information, as intended.

\begin{figure}[!b]
	 	\begin{center}
	 		\includegraphics[width=1\textwidth]{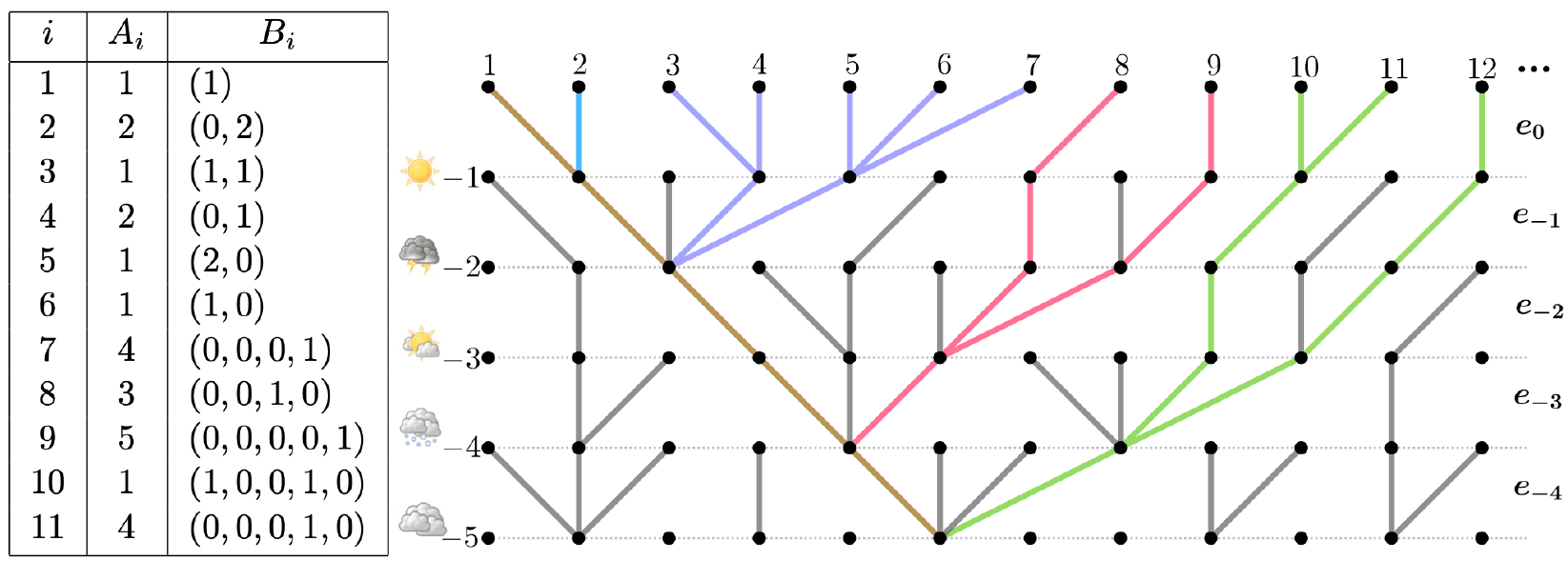}
	\caption{A Galton-Watson tree in a varying environment and the processes $(A_i, i\geq 1)$ and $(B_i, i\geq 1)$. We utilize colors to represent different subtrees rooted in the spine of individual $1$, whose descendants remain alive at the present generation. The length of vector $B_i$ corresponds to the height of the subtree attached to the first spine containing individual $i+1$.}
	\label{fig: figure-label}
	 	\end{center}
\end{figure}

We employ the same planar embedding together with the process $(D_i, i\geq 0)$ defined in \cite{amaurylea2013}. Additionally, we introduce a vector-valued process $(B_i, i\geq 0)$,  which is derived from restricting $D_i$ to some of its initial entries. Specifically, $B_i$ has length $l_i$ corresponding to the coalescent time between individuals $1,2,\dots, i+1$; with entries defined as $B_i(n)=D_i(n)$ for $1\leq n\leq l_i$. By construction, $A_i$ is the first non-zero entry of $B_i$, and $l_i=l_{i-1}\vee A_i$.
 See Figure \ref{fig: figure-label}.

If the environment remains constant, Lambert and Popovic \cite{amaurylea2013} demonstrated that the process $(A_i, i\geq 1)$ is Markov only under the condition that the offspring distribution is linear fractional, and they derived its distribution. In the case where the varying environment is linear fractional, Proposition \ref{pr:Prop Linear} demonstrates that $(A_i, i\geq 1)$ remains Markov. Furthermore, it is a sequence of independent and identically distributed random variables, and we obtain its distribution.

The remainder of the paper is structured as follows: In Section \ref{sec: preliminares}, we present Galton-Watson processes in a varying environment and discuss some related properties. Specifically, for the associated tree, we have developed the concept of stopping lines, analogous to stopping times, to establish a general branching property. The definitions of $(A_i, i\geq 1)$ and $(B_i, i \geq 1)$ are provided in Section \ref{sec: main}. Within this section, we introduce our main theorem: $(B_i, i \geq 1)$ is a Markov process that contains finite information for every $i$ and allows one to reconstruct $(A_i, i\geq 1)$. Additionally, if the offspring distributions are linear fractional, we demonstrate that the process $(A_i, i\geq 1)$ is Markov and explicitly provide its distribution. Finally, Section \ref{sec: proofs} is dedicated to the proofs.

\section{Galton-Watson processes in a varying environment}\label{sec: preliminares}

Galton-Watson processes in a varying environment model the development of the size of a population, where individuals reproduce independently, sharing a common offspring distribution that may vary across generations.

 In this context, a \textit{varying environment} is defined as a sequence $\mathcal{E}=(e_n,n\geq 1)$ of probability measures on $\N_0:=\{0,1,\dots\}$.  Let $Z_n$ represent the population size at generation $n$. The process $(Z_n,n\geq 0)$ is a Markov chain recursively defined as 
\[
Z_0=1\qquad \mbox{and} \qquad Z_{n+1}= \sum_{i=1}^{Z_n} \xi_i^{(n)}, \quad n\geq 0,
\] 
where $\xi_i^{(n)}$ denotes the number of children of the $i$-th individual living at generation $n$, with a distribution given by $e_{n+1}$. We assume that the variables  $(\xi^{(n)}_i,i\geq 1, n\geq 1)$ are all independent. The process $ \{ (Z_n, n\geq 0); \mathbb{P} \}$ is called a \textit{Galton Watson process in a varying environment} $\mathcal{E}$, abbreviated as GWVE.

Let $f_n$ be the generating function of $e_n$. For each $0\leq m< n$ and $s\in[0,1]$, we define 
\begin{equation}\label{eq:fmn}
f_{m,n}(s):=f_{m+1}\circ \cdots \circ f_n(s), \qquad \mbox{and}\qquad \qquad f_{n,n}(s):=s,
\end{equation}
where $\circ$ denotes the composition. Note that for every $s\in[0,1]$ and $0\leq m<n$,
\begin{equation}\label{eq:fgpmn}
f^\prime_{m,n}(s)=\prod_{\ell =m+1}^n f^\prime_\ell (f_{\ell,n}(s)),\qquad \mbox{and}\qquad f^\prime_{n,n}(s)=1.
\end{equation}

Let $ k \geq 0$. According to Kersting and Vatutin \cite{gotzvatutin}, we have
\begin{equation*}
\mathbb{E}(s^{Z_{k+n}} \mid Z_k=1)=f_{k, k+n}(s), \qquad 0\leq s\leq1, \quad n\geq 0.
\end{equation*}
In particular, $(Z_{k+n}, n\geq 0)$  conditionally on $\{ Z_k=1\}$, has the same law as a Galton Watson process in a varying environment $\,\mathcal{E}^{(k)}:=(e_{k+n},n\geq 1)$. To simplified the notation, we denote the law of the process with the shifted environment $\mathcal{E}^{(k)}$ by $\{ (Z_n, n\geq 0); \mathbb{P}^{(k)}\}$.

A \textit{Galton Watson tree in a varying environment} $\mathcal{E}$, abbreviated as GWVE tree, is the genealogical tree associated with $ \{ (Z_n, n\geq 0); \mathbb{P} \}$, starting with one individual. It can be viewed as a planar rooted tree $T$, with edges connecting parents to children, and the root representing the initial individual.
Any individual that appears in the tree may be labeled through its ancestry
using the Ulam-Harris notation. For example, 
consider an individual $x=(x_1,\dots ,x_{n-1},x_n)\in T$. This individual is the $x_n$-th child of the $x_{n-1}$-th child, and so forth, up to the $x_1$-th child of the root, denoted by $\emptyset$. This notation allows us to determine that the length or generation in which $x$ resides is $|x|=n$. For  two individuals $x=(x_1,\dots ,x_n)$ and $y=(y_1,\dots ,y_m)$, we denote their concatenation by $xy=(x_1,\dots ,x_n,y_1,\dots ,y_m)$. We  establish the convention that $x\emptyset=\emptyset x=x$. To define ancestry relationships within the tree, let $x,y\in T$. We say that $y$ is an ancestor of $x$ (or $x$ is a descendant of $y$) if  there exists  $w\in \{\emptyset\}\bigcup\bigcup_{n\geq 1} \mathbb{N}^n$ such that $x=yw$. This relationship is typically denoted by $y\preceq x $ (or $x\succeq y$). Observe that  $x\preceq x$, for any $x\in T$.
  
For $x\in T$, let $T^{(x )}$ denote the subtree with root in $x$. Note that $T^{(x)}$ represents the genealogical  tree associated with the process $ \{ (Z_n, n\geq 0); \mathbb{P}^{(\mid x \mid )}\}$ starting with one individual. By the branching property, it is well-known that, conditioned on the number of individuals in a given generation, the subtrees rooted at those individuals are independent. In Proposition \ref{pr:lineaparo}, we establish a stronger property. To achieve this, we introduce some concepts previously discussed in \cite{Chauvin} and later employed  by  \cite{Andreas2000} and  \cite{Jagers}.
  
 A {\it{line}} $L$ is defined as a family of vertices in $T$, such that each branch started at the root contains at most one vertex from $L$.   Consequently, for any two distinct vertices in the line, neither one is a descendant of the other. Associated with a line $L$ is its $\sigma$-algebra denoted by  $\mathcal{F}_L:= \sigma(\Omega_x: x\notin \mathfrak{D}_L)$, where
$\mathfrak{D}_L:=\{y\in T: \exists x\in L \mbox{ with } x\preceq y\}$; and for each vertex $x\in T$, $\Omega_x$ represents the information  concerning  the number of children of $x$ along with the edges connecting $x$ to its children. We also define $\mathcal{G}_x= \sigma( \Omega_y : y\preceq x)$  for any vertex $x\in T$. A {\it{stopping line}} $L$ is defined as a random line such that the event $\{x\in L  \}\in \mathcal{G}_x$, for every $x\in T$. A simple example of a stopping line is given by $L=\{ x \}$, where $x$ belongs to the set $T$. The stopping line $L^{(n)}:=\{x\in T: |x|=n\}$, with $n\in \mathbb{N}$, is commonly employed to establish the branching property. Consider an element $x\in T$ and its ancestors, denoted by $\mathfrak{a}_x:=\{y\in T: y\preceq x\}$. In the upcoming sections, we will employ a specific stopping line defined as the set containing the children of ancestors that are not ancestors themselves. This can be expressed as $L_x:=\{z\in T\setminus\mathfrak{a}_x: \exists\ y\in \mathfrak{a}_x \mbox{ with } y\preceq z \}$. Refer to the right-hand side of Figure \ref{fig:2} for an illustration, where the red dots represent the ancestors of $i$, and the stopping line is formed by the green and purple dots.

\begin{proposition} \label{pr:lineaparo}
Let $L$ be a stopping line for a GWVE tree associated with the process $ \{ (Z_n, n\geq 0); \mathbb{P} \}$. Conditioned on $\mathcal{F}_L$, the subtrees $\{ T^{(x)}: x\in L \}$ are independent. Moreover, for each $x\in L$, the tree $T^{(x)}$ conditioned on $\mathcal{F}_L$ has the same law as the genealogical tree associated with $\{ (Z_n, n\geq 0); \mathbb{P}^{(\mid x \mid )}\}$. In other words, we have
\[
\mathbb{E}\left( \prod_{x\in L} \varphi_x (T^{(x)} ) \Big|  \mathcal{F}_L \right) =  \prod_{x\in L}\mathbb{E}^{(\mid x \mid )}\left(  \varphi_x (T )   \right),
\] 
where $(\varphi_x, x\in L)$ are non-negative measurable functions on the space of planar rooted trees.
\end{proposition}
\begin{proof}
By the Monotone Class Theorem for functions, it is enough to work with indicator functions $\varphi_x= \mathbb{1}_{A_x}$ for some measurable subsets $A_x$. Let $L=\{x_1,\dots, x_n \}$ be a stopping line.  Then,  $x_i$ is neither an ancestor nor a descendant of $x_1$, for $i\geq 2$. Therefore, $T^{(x_i)}\in\mathcal{F}_{x_1}$  for  $i\geq 2$. Since $\mathcal{F}_L\subset \mathcal{F}_{x_1}$,  by the tower property  
\begin{align*}
\mathbb{E}\left( \prod_{i=1 }^n  \mathbb{1}_{A_{i}} (T^{(x_i)} )  \Big|   \mathcal{F}_L \right)
& = \mathbb{E}\left(  \mathbb{E}\left(  \prod_{i=1 }^n \mathbb{1}_{A_{i}} (T^{(x_i)} ) \Big|  \mathcal{F}_{x_1}   \right)  \Big|   \mathcal{F}_L \right)\\
& = \mathbb{E}\left(    \prod_{i=2}^n \mathbb{1}_{A_{i}} (T^{(x_i)} ) \mathbb{E}\left(   \mathbb{1}_{A_{1}} (T^{(x_1)} ) \Big|  \mathcal{F}_{x_1}  \right)  \Big|   \mathcal{F}_L \right)\\
& = \mathbb{E}^{(\mid x_1 \mid )} \left( \mathbb{1}_{A_{1}} (T^{(x_1)}  ) \right) \mathbb{E}\left(    \prod_{i=2}^n \mathbb{1}_{A_{i}} (T^{(x_i)} )   \Big|   \mathcal{F}_L \right).
\end{align*}  
In the last equality, we used the fact that $T^{(x_1)}$ is the genealogical  tree associated with  $ \{ (Z_n, n\geq 0); \mathbb{P}^{(\mid x_1 \mid )}\}$ and the branching property. Repeating the argument, we conclude the finite case. 
The general result follows from considerations of  a sequence of  finite stopping lines $L_n\uparrow L$, where $L_n$ is  monotonically increasing towards $L$, and utilizing the monotone convergence  and backwards martingale convergence theorems  applied to  the conditional  expectations given $\mathcal{F}_{L_n}$ as $\mathcal{F}_{L_n} \downarrow \mathcal{F}_{L}$.
\end{proof}

We are also interested in the law of  the children of the root  conditioned to have alive descendants at a fixed generation. We note that the survival probability up to generation $n$ is given by
\begin{equation}\label{eq:survival}
\mathbb{P}(Z_n> 0)=1 - f_{0,n}(0).
\end{equation}
For $Z_0=1$, we define $\zeta_n$ as the number of individuals at generation one having alive descendants at generation $n$. In particular, $\zeta_n$ follows the distribution
\begin{equation}\label{eq:lawzetan}
\sum_{i=1}^{\xi^{(0)}} \varepsilon_i ^{(n)},
\end{equation}
where $\xi^{(0)}$ has distribution $e_1$ and  $(\varepsilon_i^{(n)}, i\geq 1)$ is a sequence of independent Bernoulli random variables with parameter $1 - f_{1,n}(0)$, independent of $\xi^{(0)}$. Additionally, we define $\eta_n$ as $\zeta_{n}-1$ conditioned on $\{\zeta_n \geq 1\}$. 
Note that the events $\{ Z_n>0\}$ and  $\{ \zeta_n\geq1 \}$ are equivalent. Moreover, conditioned on $\{ Z_n>0\}$, there exists {\it{an}} individual at generation $1$ with alive descendants at generation $n$. Therefore, the random variable $\eta_n$ can be interpreted as the number of individuals at generation one, having alive descendants at generation $n$ which are different from {\it{that}} individual. Note that for every $k\geq 0$,
\begin{equation}\label{eq:lawetan}
\mathbb{P}(\eta_n=k) = \mathbb{P}( \zeta_n=k+1 \mid \zeta_n \geq 1 ) 
= \frac{( 1- f_{1,n}(0) )^{k+1}f_1^{(k+1)} (f_{1,n}(0))}{(k+1)! (1 - f_{0,n} ( 0 ))},
\end{equation}
where $f^{(k+1)}_1$ denotes the $(k+1)$-th derivative of $f_1$. Indeed, we use 
\[
\mathbb{E}\left(\xi^{(0)}(\xi^{(0)}-1)\cdots (\xi^{(0)}-k)s^{\xi^{(0)}-k+1}\right)=f^{(k)}_1(s)
\] 
to obtain 
\begin{align*}
\mathbb{P}(\zeta_n=k+1)
& = \sum_{i=1}^\infty \mathbb{P}( \zeta_{n}=k+1 \mid  \xi^{(0)}=i) \mathbb{P}(\xi^{(0)}=i) \\
& = \sum_{i=1}^\infty \binom{i}{k+1} (1 - f_{1,n}(0))^{k+1} (f_{1,n}(0))^{i-k-1}  \mathbb{P}(\xi^{(0)}=i) \\
& =\frac{(1 - f_{1,n}(0))^{k+1} }{ (k+1)! }  f_1^{(k+1)}( f_{1,n}(0) ).
\end{align*}
Similarly, 
\begin{align*}
\mathbb{P}(  \zeta_n=0) 
& = \sum_{i=0}^\infty (  f_{1,n}(0) )^i q_1(i) 
= f_1( f_{1,n}(0) ) =f_{0,n} ( 0 ).
\end{align*}

We complete the section with an example in which we compute the aforementioned  distribution when the offspring distributions are linear fractional.

\begin{example}\label{Ex linear1} An offspring distribution $\xi$ that satisfies
	\[
	\mathbb{P}(\xi =0)=1-r \qquad \mbox{and} \qquad \mathbb{P}(\xi =k) =rpq^{k-1}, \quad k\neq 0,
	\] 
	where $0<q<1$, $p=1-q$ and $0\leq r \leq 1$, is called a \textit{linear fractional} distribution with parameters $(r,p)$. Special cases include $r=1$, yielding a geometric distribution $G$, and $r=0$, yielding a Dirac measure $\delta_0$. If $0<r<1$, the distribution is a mixture of both, i.e.  $\xi=rG+(1-r)\delta_0$. In this case, its generating function is given by
	\[
	f(s)=1-r \frac{1-s}{1-qs}, \qquad s\in [0,1].
	\] 
	Thanks to the identity
	\begin{equation*}
	\frac{1}{1-f(s)}=\frac{1}{1-f'(1)(1-s)} + \frac{ f''(1)}{2f'(1)^2}, \qquad s\in[0,1],
	\end{equation*}
	we can see that a linear fractional distribution is characterized by  its mean $f'(1)=r/p$ and its normalized second factorial moment $f''(1)/f'(1)^2=2q/r$. 
	
	We say that a \textit{varying environment} $\mathcal{E}=(e_n,n \geq 1)$ is \textit{linear fractional} if and only if every $f_n$ is linear fractional with parameters $\{(r_n,p_n): n\geq 1 \}$.  According with \cite[Chapter 1]{gotzvatutin}, the generating function $f_{m,n}$ is again linear fractional with mean
	\[
	f_{n,n}'(1)=1,\qquad	f_{m,n}'(1)=f_{m+1}'(1)\cdots f_n'(1) = \frac{r_{m+1} \cdots r_n}{p_{m+1} \cdots p_n}, \qquad m<n,
	\] 
	and normalized second factorial moment \  $f_{n,n}''(1)/f_{n,n}'(1)^2=0$; and for $m<n$,
	\[	
	\frac{f_{m,n}''(1)}{f_{m,n}'(1)^2}
	=\frac{f_{m+1}''(1)}{f_{m+1}'(1)^2}+\underset{k=m+2}{\overset{n}{\sum}} \frac{1}{f_{m+1}'(1)\cdots f_{k-1}'(1)}\frac{f_k''(1)}{f_k'(1)^2} =\frac{2q_{m+1}}{r_{m+1}}+\hspace{-.2cm}\underset{k=m+2}{\overset{n}{\sum}} \frac{ 2p_{m+1}  \cdots p_{k-1} q_k}{r_{m+1}  \cdots r_{k-1}r_k}.
	\] 
	Consider a Galton-Watson process $(Z_n:n\geq 0)$ with the aforementioned  environment. Then, $Z_n$ is linear fractional with generating function $f_{0,n}$. Furthermore, by using \eqref{eq:lawetan} and
	\[
	f_1^{(k+1)}(s) = \frac{(k+1)! r_1(1-q_1)q_1^{k} }{ (1-q_1 s)^{k+2} },
	\] 
	we can prove that  $\eta_{n}$  follows a geometric distribution. More precisely,  
	\begin{equation*}
	\begin{split}
	\mathbb{P}(\eta_n=k) =
	& \left( \frac{1-q_1}{1-q_1  f_{1,n}(0)} \right)  \left( \frac{q_1 (1 - f_{1,n}(0) )}{1-q_1  f_{1,n}(0)} \right)^k,\qquad k\geq 0.
	\end{split}
	\end{equation*}
	In other words, 	$\mathbb{P}(\eta_1=k)=(1-q_1)q_1^k,$ for $ k\geq 0$; and for $n\geq 2$ 
	\begin{equation*}
	\begin{split}
	\mathbb{P}(\eta_n=k) =
	&
	\left( \frac{ 1 + \underset{i=2}{\overset{n}{\sum}}
s_{i,n} }{  1 + \underset{i=1}{\overset{n}{\sum}}
		 s_{i,n}}  \right)
	\left( \frac{  s_{1,n}  }{  1 + \underset{i=1}{\overset{n}{\sum}}
		s_{i,n} }  \right) ^k, \qquad k\geq 0,
	\end{split}
	\end{equation*}
	where $s_{n,n}=\tfrac{1-p_n}{p_n}$ and $s_{i,n}=\tfrac { (1-p_i) }{p_{i}}\tfrac{r_{i+1}  \cdots r_n }{  p_{i+1}\cdots p_n}$ for $\,i<n$.
\end{example}

\bigskip
In the next section, we will analyze the backward genealogy of a random population when the forward in time dynamics are produced by a Galton-Watson process in a varying environment.   The offspring distributions will be denoted by $\mathcal{E}=(e_m, m\in \Z_{-})$. The generating functions $f_{m,n}$ in \eqref{eq:fmn} can be extended to include $m\leq n \leq 0$.

\section{Main results}\label{sec: main}

Consider an arbitrarily large population at the present time, originating from an unspecified time in the distant past, where individuals within the same generation reproduce independently forward in time, sharing the same offspring distribution that may vary across generations. We denote by $e_{m+1}$ the offspring distribution of individuals at generation $m$. It is essential to note that for each individual in the past, the forward-in-time process starting from that individual generates a branching tree in a varying environment.

To analyze the backward genealogy of the present population, we adopt a specific embedding of a branching tree, which is infinite in the number of generations and
individuals at any generation. This embedding serves as the basis for defining the coalescent point process, initially proposed by Lambert and Popovic \cite{amaurylea2013} for a Galton Watson process, and later extended by Popovic and Rivas \cite{LeaRivas2014} to multitype Galton Watson processes. In this planar embedding, individuals are situated at points $(m,i)$ on a discrete lattice $\mathbb{Z}\times \mathbb{N}$, where the first coordinate $m$ represents the generation, and the second coordinate $i$ represents the individual's position in the planar embedding layout from left to right. Given our focus on backward genealogies, our analysis is confined to $\mathbb{Z}_{-}$, where the present generation is denoted by $m=0$. Each vertex in the lattice is connected to its offspring, represented by vertices in
the level above, in a manner that avoids empty spaces and intersections between lineages. Specifically, let $\xi^{(m)}_i$ denote the number of offspring of individual $(m,i)$.  We define the population's genealogy such that individual $(m,i)$ has mother $(m-1,j)$ if and only if
\[
\sum_{k=1}^{j-1} \xi_{k}^{(m-1)} < i \leq \sum_{k=1}^j \xi_{k}^{(m-1)}.
\] 
We suppose that the variables $\{\xi^{(m)}_i, i\geq 1, m\in \Z_{-}\}$ are independent, and  the sequence $\{\xi^{(m)}_i, i\geq 1\}$  has common distribution $e_{m+1}$ for every $m\in \Z_{-}$.  In the remainder of the paper, the varying environment  is  denoted by $\mathcal{E}=(e_m, m\in \mathbb{Z}_{-})$. 

We now proceed, following the approach and notation outlined in \cite{amaurylea2013}, by introducing the elements necessary for recovering the genealogy. For each $n \in \mathbb{N}$, let $\mathfrak{a}_i(n)$ denote the index of the ancestor of individual $i$ in generation $-n$. The coalescent time $C_{i,j}$ of individuals $i$ and $j$ is denoted by
\[
C_{i,j}:=\min \{ n\geq 1: \mathfrak{a}_i(n)=\mathfrak{a}_j(n)  \}, \qquad i,j\in \N,
\] 
with the understanding that $\min \emptyset = \infty$. We define 
\[
A_i:=C_{i,i+1}, \qquad i\in\mathbb{N}.
\] 
By construction, $C_{i,j}=\max\{ A_i,A_{i+1},\dots, A_{j-1}\}$, for any $i<j$. Therefore, $(A_i, i \geq 1)$ contains all the genealogical information of the current population. The sequence $\mathbf{A}:=(A_i, i \geq 1)$  is called the \textit{coalescent point process in varying  environment}. 
 This process was first defined in \cite{amaurylea2013}, when the forward-time dynamics were produced by a branching process in a constant environment. 

 The distribution of $\mathbf{A}$ is not easy to determine, and in general, it is not a Markov process, except for some special cases.  Following the approach of \cite{amaurylea2013}, let us define an auxiliary Markov process that characterizes the genealogy. For any fixed individual $i$, we follow its ancestral line (referred to as the $i$-th spine) and consider the subtrees with roots in this spine. Note that these roots are labeled by $\{(-n,\mathfrak{a}_i(n)), n\in\mathbb{N} \}$. At every subtree, we count the number of daughters of the root at the right-hand side of the spine, whose descendants are alive at the present generation. To be precise,  let 
\[
\mathcal{D}_i(n):=\{ \text{daughters of } (-n,\mathfrak{a}_i(n)) \text{ with descendants in } \{ (0,j): j\geq i  \}  \}, \qquad n\geq 1, \ i\geq 1.
\] 
Define
\[
D_i(n) = \# \mathcal{D}_i(n)-1, \qquad i,n\geq 1.
\] 
For every $i\geq 1$, we define the sequence $D_i:=(D_i(n),n\geq 1)$. We set $D_0$ as the null sequence and $A_0:=\infty$. It follows from the monotone planar embedding that
\begin{equation} \label{eq:AdeD}
A_i=\min\{n\geq 1: D_i(n)\neq 0\}, \qquad i\geq 0.
\end{equation}
We will establish in Proposition \ref{pr:ThD} that  the  sequence-valued process ${\mathbf{D}}:=(D_i, i\geq 0)$ possesses the Markov property. To provide the transitions probabilities,  we begin by recalling that for each $(m,i)\in \mathbb{Z}_{-}\times \N$,  the forward-in-time dynamics starting at individual $(m,i)$  are produced by a branching process in a varying environment. To be precise, let  $Z^{(m,i)}(k)$ denote the number of descendants of individual $(m,i)$ at generation $m+k$. Then,  $(Z^{(m,i)}(k), k\geq 0)$ constitutes a GWVE process with environment $\mathcal{S}_{(m)}:=(e_{m+1}, e_{m+2}, \dots, e_0, 0, 0, \dots)$, starting with one individual. Consequently, for every  $m\in \Z_-$, the probability that an individual at generation $m$ has alive descendants at the present generation is $p_{m}=1-f_{m,0}(0)$, see \eqref{eq:survival}.  Moreover, if we denote by  $\zeta^{(m)}$ the number of its daughters with alive descendants at generation zero, as a consequence of \eqref{eq:lawzetan}, we have 
	\begin{equation*}
		\zeta^{(m)}\overset{\mathcal{L}}{=}\underset{i=1}{\overset{Y}{\sum}}\varepsilon_i,
	\end{equation*}
	where $Y\sim e_{m+1}$ and the variables $\varepsilon_i$ are Bernoulli with parameter $1-f_{m+1,0}(0)$, all independent. We also define $\eta^{(m)}$ as
	\begin{equation}\label{eq:etadistribucion}
		\eta^{(m)}=\zeta^{(m)}-1, \qquad \mbox{ conditioned on } \{\zeta^{(m)}>0\}.
	\end{equation}
	In particular, by \eqref{eq:lawetan}
	\begin{equation}\label{eq:etanueva}
		\mathbb{P}(\eta^{(m)}=0) = \frac{ ( 1- f_{m+1,0}(0) )f_{m+1}^\prime (f_{m+1,0}(0))}{ 1 - f_{m,0} ( 0 )}.
	\end{equation}
Now, we are ready to establish the transition probabilities of the process $(D_i,i\geq 0)$. This is an extension of \cite[Theorem 2.1]{amaurylea2013} to varying environment, and we will use similar techniques to prove it.

\begin{proposition}\label{pr:ThD}
	The sequence-valued process $(D_i,i\geq 0)$  is  a Markov chain starting at the null sequence with transition probabilities given by
	\begin{equation*}
	\Big( D_{i+1}(m) \ \mid\  D_i = (d(n))_{ n\geq 1 } \Big) \ {=}
	\left\{ 
	\begin{array}{ll}
	\eta^{(-m)}  & \text{ if } 1\leq m <A_i\\
	d(m)-1& \text{ if }m=A_i\\
	d(m) & \text{ if } A_i<m,
	\end{array}
	\right.
	\end{equation*}
	where $(\eta^{(-m)}, m\geq	 1 )$ is a sequence of independent random variables such that $\eta^{(-m)}$ is distributed as \eqref{eq:etadistribucion}, for each $m$. Moreover, the law of $A_1$ is given by
	\[
	\mathbb{P}(A_1>n)=  \prod_{i=1}^n \mathbb{P}(\eta^{(i)}=0)=\frac{ f^{\prime}_{-n,0}(0) }{1-f_{-n,0}(0)}= \mathbb{P}(Z_n=1\mid Z_n=0), \qquad 
	\]
	where $(Z_k, k\geq 0)$ is a GWVE with environment $\mathcal{S}_{(-n)} $ for every $n\geq 1$.
\end{proposition}

Considerable repetitive information exists within $(D_i,i\geq 0)$, as for every $i$, $D_i$ and $D_{i+1}$ are essentially identical infinite sequences, differing only in a finite number of entries. To address this, a point measure-valued process $(\widetilde{B}_i, i \geq 0)$ was recursively defined in \cite{amaurylea2013} for 
 a constant environment. In essence, when {$\widetilde{B}_i$ has positive measure at $n$,  its  multiplicity records the number of children of $(-n,\mathfrak{a}_i(n))$ located at the right-hand side of the $i$-th spine, with descendants among individuals $\{(0,j) : j\geq i+1\}$}. The construction proceeds as follows: $\widetilde{B}_0$ is the null measure. $\widetilde{B}_1$  has positive mass at position $A_1$, where its multiplicity records the number of children of individual $(-A_1, \mathfrak{a}_1(A_1))$ with descendants in individuals  $\{(0,j) : j\geq 2\}$. Recursively, $\widetilde{B}_{i+1}$ is updated from $\widetilde{B}_i$ by reducing by one the mass at position $A_i$ (as the $(i+1)$-th spine is part of the children at the right hand side of the $i$-th spine) and possibly by adding a new mass at position $A_{i+1}$ with the respective multiplicity. Formally, for any finite point measure $b=\sum_{n\geq 1}b(n)\delta_n$, we define the minimum of its support as
\[
\mathfrak{s}( b ) := \min\{n\geq 1: b(n)\neq 0\}.
\] 
Additionally, we define  $b^\ast=b-\delta_{\mathfrak{s}( b )},$ with the convention that $\mathfrak{s}( b )=\infty$ and  $b^\ast=b$  if $b$ is the null measure. Then, the process $\widetilde{B}_i$ is recursively defined as follows  
\begin{equation*}
\widetilde{B}_{i+1} :=
\left\{ 
\begin{array}{ll}
\widetilde{B}_i^\ast+D_{i+1}(A_{i+1})\delta_{A_{i+1}} & \ \text{ if } A_{i+1}\neq \mathfrak{s}( \widetilde{B}_i ) \text{ and }   A_{i+1}< \mathfrak{s}( \widetilde{B}_i^\ast ),\\
\widetilde{B}_i^\ast& \ \text{ otherwise.}
\end{array}
\right.
\end{equation*}

Here, we present a counterexample to illustrate that $(\widetilde{B}_i, i \geq 0)$ may not constitute a Markov process, thereby challenging the validity of \cite[Theorem 2.2]{amaurylea2013}.
\begin{example}\label{counterexample}
	Consider the process $(\widetilde{B}_i,i\geq 1)$ associated to Figure \ref{fig: figure-label}: 
	\begin{center}
	\begin{tabular}{l l l c c c }
	$\widetilde{B}_1=\delta_1$, 
	& $\widetilde{B}_2=2\delta_2$, &$\widetilde{B}_3=\delta_1+\delta_2$, 
	& $\widetilde{B}_4=\delta_2$,
	&$\widetilde{B}_5=2\delta_1$,
	& $\widetilde{B}_6=\delta_1$,\\
	$\widetilde{B}_7=\delta_4$,
	& $\widetilde{B}_8=\delta_3$, &$\widetilde{B}_9=\delta_5$, &$\widetilde{B}_{10}=\delta_1$ 
	& and &$\widetilde{B}_{11}=\delta_4$. 
	\end{tabular}
\end{center}
We show that this process  \textit{cannot be  Markov} by exhibiting that the value of $\widetilde{B}_{11}$ depends on both  $\widetilde{B}_{9}$ and $\widetilde{B}_{10}$. Through the recursive construction, we observe that $\widetilde{B}_{10}=\delta_1$ implies that $\widetilde{B}_{11}=\widetilde{B}_{11}(N)\delta_N$ for some  $N\geq 2$ with $\widetilde{B}_{11}(N)>0$. Furthermore, from $\widetilde{B}_9=\delta_5$ we know that individual $(-5,\mathfrak{a}_9(5))$  has a \textit{unique} daughter $u$, at the right hand side of the $9$-th spine, whose descendants are alive at the present generation. Since $A_{9}=5$ and $A_{10}=1$, we have $(-5,\mathfrak{a}_9(5))=(-5,\mathfrak{a}_{10}(5))=(-5,\mathfrak{a}_{11}(5))$  and $u$ is part of both the $10$-th  and $11$-th spine.  Therefore, $\widetilde{B}_{11}(5)=0$ and $N$ cannot be $1$ or $5$.  This implies that  the value of $\widetilde{B}_{11}$  depends not only on $\widetilde{B}_{10}$ but also on  $\widetilde{B}_{9}$. \hfill \qed
\end{example}

Next, we focus on constructing our process $\mathbf{B}=(B_i, i\geq 0)$. For each $i$, we  restrict $D_i$ to the first $C_{1,i+1}$ entries, where we recall that $C_{1,i+1}$ denotes the first generation in which an individual in the $i$-th spine is an ancestor of the individuals in the present generation indexed from $1$ to $i+1$. Equation \eqref{eq:AdeD}, recursively defines the length the length $C_{1,i+1}$ in terms of $D_i$ as
\begin{equation} \label{eq:defli}
C_{1,i+1}= C_{1,i}\vee A_i =C_{1,i} \vee \min\{n\geq 1: D_i(n)\neq 0\}.
\end{equation}
We define the process $\mathbf{B}=(B_i, i\geq 0)$ as $B_0=\emptyset$ and $B_i:=(D_i(1), D_i(2), \dots, D_i(C_{1,i+1}))$ for $i\geq 1$. Henceforth, we refer to the process $(B_i, i\geq 0)$ as the \textit{coalescent point process with multiplicities in a varying environment}. See Figure \ref{fig: figure-label} for an illustration. 

The main result of this work establishes that the process $\mathbf{B}$ is Markov with a state space $\mathcal{V}$, defined as the set of all the vectors with entries in $\mathbb{N}$, i.e.
	\[
	\mathcal V = \bigcup_{m\in\mathbb{N}} \mathbb{N}^m
	\] 
	with the convention  that $\mathbb{N}^0=\emptyset$. For any $b=(b(1), b(2), \dots ,b(m)) \in \mathcal{V}$, its length is denoted by $l(b)=m$ with the convention   $l(\emptyset)=0$. The function $\mathfrak{s}(b)$ indicates the first non-zero coordinate of $b$, with the convention $\mathfrak{s}( b ) =\infty$, if $b$ is a null vector or the empty set. For any $b=(b(1), b(2), \dots ,b(m)) \in \mathcal{V}$, we define the vector $b^*=(b^*(1), b^*(2), \dots ,b^*(m))$ by
	\begin{equation*}
	b^*(j)= b(j) - \mathbb{1}_{ \{j=\mathfrak{s}(b) \}  },\qquad j\leq m,
	\end{equation*}
	with the convention $b^*=b$, if $b$ is a null vector or the empty set. It is worth noting that in our definition of  $(B_i, i\geq0)$, we use indicator functions  instead of Dirac measures. This choice arises from the fact that our process $(B_i, i\geq0)$ is vector-valued, contrasting with the point-measure-valued nature of the process $(\widetilde{B}_i, i\geq0)$ in \cite{amaurylea2013}. Furthermore, it is important to highlight that in the construction of $(\widetilde{B}_i, i\geq0)$, only the mass at position $A_i$ could be added, whereas in our process, masses above level $C_{1,i}$ could be added.

\begin{theorem}\label{th:ThB}
The vector-valued process $\mathbf{B}=(B_i,i\geq 0)$ is a Markov chain starting at $B_0=\emptyset$. Conditioned on the event $\,\{ B_i=(b(1),\dots, b(\ell))\}$, the law of the vector $B_{i+1}$ is determined by the following transition probabilities
	\begin{equation*}
	B_{i+1}(m):  =
	\left\{ 
	\begin{array}{ll}
	\eta^{(-m)} & \text{ if }\, 1\leq m < A_i\ \text{  or  }\  \ell< m\leq  l(B_{i+1}) \\
	b(m)-1  & \text{ if }\, m= A_i\\
	b(m)& \text{ if }\, A_i < m \leq \ell  ,
	\end{array}
	\right.
	\end{equation*}
	where $(\eta^{(-m)}, m\geq 1 )$ is a sequence of independent random variables such that for each $m$, the variable $\eta^{(-m)}$ is distributed as \eqref{eq:etadistribucion}. The length of $B_{i+1}$ satisfies
	\begin{equation}\label{eq:LMarkov}
l(B_{i+1}) = \ell \vee \left(\mathbb{1}_{ \{ \mathfrak{s}(B_i^*)=\infty\}}\min\{k \in \{1,2,\dots , A_i -1\}\cup \{ \ell +1,\ell +2, \dots \} : \eta^{(-k)} \neq 0\}\right). 
	\end{equation}
Moreover,  $A_i$ is a functional of $B_i$, i.e. if $B_i=(b(1),\dots, b(\ell))$, then
\begin{equation} \label{eq:AdeB}
A_i=\min\{ 1 \leq n\leq \ell : b(n)\neq 0\}, \qquad i\geq 1.
\end{equation}
\end{theorem}
	
To conclude this section, we establish that if the offspring follows a  linear fractional distribution, the coalescent point process $(A_i,i\geq 0)$ is Markov, and we derive its distribution. This observation was previously noted in \cite[Proposition 5.1]{amaurylea2013} for a constant  environment, and in \cite{rannala1997genealogy} with an alternative formulation.

\begin{proposition}\label{pr:Prop Linear} Assume that the environment $\mathcal{E}$ is linear fractional with parameters  $\{ (r_m, p_m): m\in \mathbb{Z}_-\}$. Then, the variables $(A_i,i\geq 1)$  are independent with common distribution given by
	\begin{equation*}
	\mathbb{P}(A_1>n)= \left( 1 + \underset{i=-n+1}{\overset{0}{\sum}}
	s_{i} \right)^{-1}, \qquad n\geq 1,
	\end{equation*}
	where $s_{0}=\tfrac{1-p_0}{p_0}$ and \ $s_{i}=\frac {  1- p_i}{p_{i}}\tfrac{ r_{i+1}  \cdots r_0 }{  p_{i+1}\cdots p_0}  $, \  for \, $i<0$.
\end{proposition}

\section{Proofs }\label{sec: proofs}
In this section, we present the proofs of  Proposition \ref{pr:ThD}, Theorem \ref{th:ThB} and Proposition \ref{pr:Prop Linear}. It is worth noting that Propositions  \ref{pr:ThD} and  \ref{pr:Prop Linear} extend  \cite[Theorem 2.1 and Proposition 5.1]{amaurylea2013} to a varying environment. Throughout the proofs, we will use some of the techniques employed in their work. Additionally, we will use a connection between the planar rooted trees introduced in Section \ref{sec: preliminares} and the monotone planar embedding described in Section \ref{sec: main}.

For a  fixed $(m,i)\in \mathbb{Z}_{-}\times \N$, recall that  $(Z^{(m,i)}(k),k\geq 0)$  is a GWVE with environment $\mathcal{S}_{(m)}=(e_{m+1}, e_{m+2}, \dots, e_0, 0, 0, \dots)$, starting with one individual. Here, $Z^{(m,i)}(k)$ denotes the number of descendants of individual $(m,i)$ at generation $m+k$.  According to Section \ref{sec: preliminares}, its genealogical tree can be viewed as a planar rooted tree denoted as $\mathcal{T}^{(m,i )}$. To simplify, we use $ \{ (Z_k, k\geq 0); \mathbb{Q}_{(m)} \}$ to denote a GWVE with environment $\mathcal{S}_{(m)}$.

In accordance with Section \ref{sec: preliminares}, for every $x\in \mathcal{T}^{(m,i )}$, the subtree with root at $x$, denoted as $\mathcal{T}^{(x )}$, shares the same distribution as the tree associated with  $ \{ (Z_k,k\geq 0); (\mathbb{Q}_{(m)})^{ (|x|)} \}$, where $|x|$ denotes its length in  $\mathcal{T}^{(m,i )}$. It is worth noting that  $0\leq |x|\leq -m$. Consequently, by shifting the environment, $\mathcal{T}^{(x )}$ shares the same law as the tree associated with  $ \{ (Z_k,k\geq 0); \mathbb{Q}_{(m+|x|)} \}$.

 Let $L$ be a stopping line for $ \mathcal{T}^{(m,i )}$. We recall that $\mathcal{F}_L= \sigma(\Omega_x: x\notin \mathfrak{D}_L)$, where for each $x\in \mathcal{T}^{(m,i )}$,  $\Omega_x$ is the information about the number of children of $x$ together with the edges connecting $x$ with its children; and $\mathfrak{D}_L$ is the set of descendants from individuals in the line. Remember that $L\subset \mathfrak{D}_L$.
 
 Based on Proposition  \ref{pr:lineaparo},  conditioned on $\mathcal{F}_L$, the subtrees $\{ \mathcal{T}^{(x)}: x\in L \}$ are independent. Furthermore,
\begin{equation}\label{eq: propapl}
\mathbb{P}\left( \bigcap_{x\in L} \{\mathcal{T}^{(x)}  =  {\bf t}_{x}\}  \mid  \mathcal{F}_L \right) =  \prod_{x\in L} \mathbb{Q}_{(m+|x|)}  \left(    {\bf t}_{x}  \right),
\end{equation}
where  ${\bf t}_{x} $  is a (deterministic) rooted tree, for every $x\in L$.

\subsection{Proof of Proposition \ref{pr:ThD}}	
 Since $A_i$ is the coalescent time of individual $i$ and $i+1$, the $i$-th and $(i+1)$-th spines coincide for every  generation $-m$ where $m\geq A_i$. This implies,
	\begin{equation}\label{eq frak a}
	\mathfrak{a}_i(m) \neq\mathfrak{a}_{i+1}(m)   \quad  \text{for all } m< A_i \qquad 
	\text{and} \qquad \mathfrak{a}_i(m) =\mathfrak{a}_{i+1}(m)   \quad  \text{for all } m\geq A_i.
	\end{equation}

	Recall that $D_i(m)$ counts the number of children of $(-m, \mathfrak{a}_i(m))$ at the right hand side of the $i$-th spine, whose descendants are alive at the present generation. Observe that $(-A_i +1, \mathfrak{a}_{i+1}(A_i -1)  )$ is a child of $(-A_i, \mathfrak{a}_{i}(A_i )  )  = (-A_i, \mathfrak{a}_{i+1}(A_i )  )$. Moreover, it is part of the $(i+1)$-spine and it is at the right hand side of the $i$-th spine. Then,  $D_{i+1}(A_i) = D_i(A_i) -1$.

	By definition of $A_i$, for every $m> A_i$, each daughter of $(-m, \mathfrak{a}_{i+1}(m))$ has descendants in $\{(0,j): j\geq i+1\}$ if and only if  she has descendants in $\{(0,j): j\geq i\}$.  Then, $D_i(m)=D_{i+1}(m)$ for all $m>A_i$.   
	
\begin{figure}[!b]
	 	\begin{center}
	 		\includegraphics[width=1\textwidth]{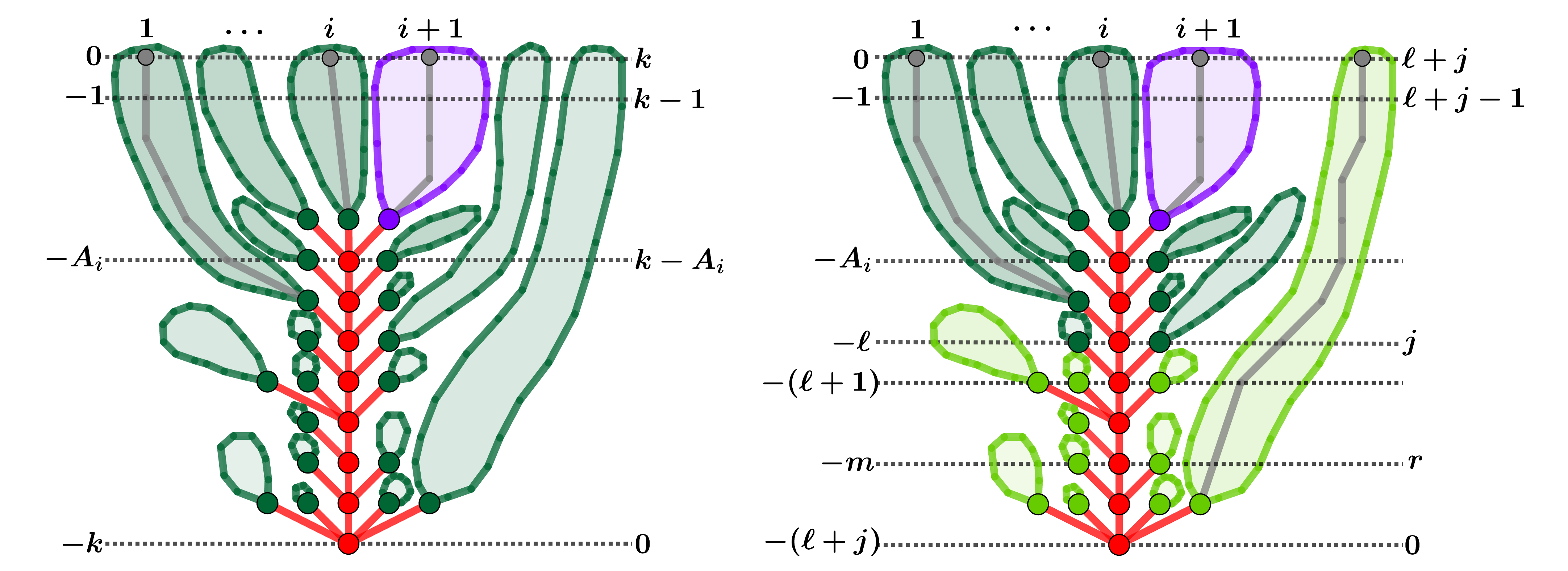}
	 		\caption{ At the left hand-side of both figures, we have the generation label as in the planar embedding $\mathbb{Z}_{-}\times\mathbb{N}$. At the right hand-side of both figures, we have the generation label as in a GWVE tree. The left tree corresponds to  $\mathcal{T}^{(-k , \mathfrak{a}_i( k) )}$.  The roots in the green and purple subtrees represent the set $L_k$. The $\sigma$-algebra $\mathcal{F}_{L_k}$ is generated by the edges and vertices in red. $R_k$ is the set of green roots.  The $\sigma$-algebra $\mathcal{G}_k$ is generated by the trees rooted in $R_k$. In the right tree, we consider $k=\ell +j$.  The set $L_{\ell+j}$ can be decomposed into the purple root, together with  
$R_{\ell +1}$ colored in dark green, and $L_{\ell+j}\setminus L_{\ell+1}$ colored in light green. }
	 		\label{fig:2}
	 		%\hspace{-.75cm}	t
	 	\end{center}
	 \end{figure}

Now, let us analyze the case $1\leq m < A_i$. Define the set
 	\[
	G_i:=\{ D_j(n)=d_j(n): n\geq 1 \, \mbox{ and }\, 1\leq j\leq i \}.
	\]
where $\{d_j(n):n\geq 1 \, \mbox{ and }\, 1\leq j\leq i \}$ are fixed integers. For every $k> A_i$, consider the tree 
$\mathcal{T}^{(-k , \mathfrak{a}_i( k) )}$, which is a GWVE tree with environment $\mathcal{S}_{(k)}$ starting from individual $(-k, \mathfrak{a}_i( k) )$. Let $L_k$ be  the stopping line defined by 
\begin{equation}\label{def:stoppingline}
L_k:=\{\mbox{daughters of } (-n,\mathfrak{a}_{i}(n)): A_i \leq  n\leq k \}\setminus \{(-n,\mathfrak{a}_{i}(n)): A_i \leq  n\leq k-1  \}.
\end{equation}
Note that
$\mathcal{F}_{L_{k}}=\sigma(\Omega_{(-n,\mathfrak{a}_{i}(n))}: A_i \leq  n\leq k )$.  Additionally, define  the set of roots
\begin{equation}\label{def:roots}
R_k:=L_k\setminus \{ (-A_i+1,\mathfrak{a}_{i+1}(A_{i}-1 ))\}.
\end{equation}
Let $\mathcal{G}_k$ denote the $\sigma$-algebra generated by the trees with roots in $R_k$. Refer to the left tree in Figure \ref{fig:2} for a representation.
According to Proposition \ref{pr:lineaparo}, conditioned on $\mathcal{F}_{L_{k}}$, the trees $\mathcal{T}^{(x)}$ with roots $x\in R_k$ and the tree rooted at $(-A_i+1,\mathfrak{a}_{i+1}(A_{i}-1)  )$ are independent. In other words, conditioned on $\mathcal{F}_{L_{k}}$, the tree $\mathcal{T}^{(-A_i+1,\mathfrak{a}_{i+1}(A_{i}-1)  )}$ and $\mathcal{G}_k$ are independent. Then, by Proposition III.3.2 in \cite{Cinlar}  
\begin{equation}\label{eq:cinlar}
\mathbb{P}\left(   \mathcal{T}^{(-A_i+1,\mathfrak{a}_{i+1}(A_{i}-1)  )} ={\bf t} \mid  \mathcal{G}_k \vee \mathcal{F}_{L_{k}}\right) = \mathbb{P}\left(   \mathcal{T}^{(-A_i+1,\mathfrak{a}_{i+1}(A_{i}-1)  )} = {\bf t} \mid \mathcal{F}_{L_{k}} \right),
\end{equation}
for every deterministic rooted tree $\bf{t}$. We observe that $\mathcal{G}_k\subset \mathcal{G}_{k+1}$ and $\mathcal{F}_{L_k}\subset \mathcal{F}_{L_{k+1}}$ for all $k$. Let $\mathcal{G}:=\lim_k \mathcal{G}_k$  and $\mathcal{F}:=\lim_k\mathcal{F}_{L_{k}}$. Note that $G_i$ is measurable with respect to $\mathcal{G} \vee \mathcal{F}$. Therefore, by the Monotone Convergence Theorem and the previous equation,
\begin{align*}
\mathbb{P}(  \mathcal{T}^{(-A_i+1,\mathfrak{a}_{i+1}(A_{i}-1)  )} = {\bf t} \mid G_i)
& = \mathbb{P}(  \mathbb{P} ( \mathcal{T}^{(-A_i+1,\mathfrak{a}_{i+1}(A_{i}-1)  )} ={\bf t} \mid \mathcal{G} \vee \mathcal{F} )  \mid G_i) \nonumber \\
&= \mathbb{P}(  \lim_k\mathbb{P} ( \mathcal{T}^{(-A_i+1,\mathfrak{a}_{i+1}(A_{i}-1)  )} ={\bf t} \mid \mathcal{G}_k \vee \mathcal{F}_{L_{k}} )  \mid G_i) \nonumber \\
& = \mathbb{P}( \lim_k  \mathbb{P} ( \mathcal{T}^{(-A_i+1,\mathfrak{a}_{i+1}(A_{i}-1)  )} = {\bf t} \mid   \mathcal{F}_{L_{k}}  )    \mid G_i).
\end{align*}
Now,  we use \eqref{eq: propapl}, with $m=-k$,   $x=(-A_i+1,\mathfrak{a}_{i+1}(A_{i}-1)  )$ and $|x|=k-A_i+1$,  to obtain
\begin{equation}\label{eq:medq}
\mathbb{P} ( \mathcal{T}^{(-A_i+1,\mathfrak{a}_{i+1}(A_{i}-1)  )} = {\bf t} \mid   \mathcal{F}_{L_{k}}  ) =\mathbb{Q}_{(-A_i+1)}  \left(    {\bf t} \right), \qquad \mbox{ for all } \quad k> A_i.
\end{equation}
Hence,
\begin{equation*}
\mathbb{P}(  \mathcal{T}^{(-A_i+1,\mathfrak{a}_{i+1}(A_{i}-1)  )} = {\bf t}  \mid G_i) =  \mathbb{P}(  \mathbb{Q}_{(-A_i+1)}  \left(    {\bf t} \right)   \mid G_i).
\end{equation*}
Conditional on $G_i$, observe that $A_i$ only depends on the values of $D_i$. Moreover, the vector $(D_{i+1}(1), D_{i+1}(2), \dots, D_{i+1}(A_{i} -1) )$ is a function of the tree $\mathcal{T}^{(-A_i+1,\mathfrak{a}_{i+1}(A_{i}-1)  )}$. Therefore, $(D_{i+1}(1), D_{i+1}(2), \dots, D_{i+1}(A_{i} -1) )$ given $G_i$ only depends on the values of $D_i$. This implies the Markov property. Now, let us determine its distribution.

Let $I(n,j)$ be the smallest index among all the descendants of individual $(-n,\mathfrak{a}_j(n))$ at the present generation, 
\[
I(n,j) :=\min\{k\leq j : \mathfrak{a}_k(n) =\mathfrak{a}_j(n)  \}, \quad j,n\geq 1.
\]
Note that  given $G_i$,   $I(-A_{i}+1, \mathfrak{a}_{i+1}(A_i-1)) =i+1$. Then,  conditioned  on $G_i$, for every $1\leq m <A_i$,
\begin{equation}\label{eq:same}
\begin{split}
\mathcal{D}_{i+1}(m) 
& =  \{ \text{daughters of } (-m,\mathfrak{a}_{i+1}(m)) \text{ with descendants in } \{ (0,j): j\geq i+1  \}  \} \\
& = \{ \text{daughters of } (-m,\mathfrak{a}_{i+1}(m)) \text{ with descendants in } \{ (0,j): j\geq 1  \}  \}.
\end{split}
\end{equation}
It follows that for every $1\leq m <A_i$, the variable $D_{i+1}(m)$ conditioned on $G_i$  has the same  distribution as $\eta^{(-m)}$  given in \eqref{eq:etadistribucion}. Moreover, $D_{i+1}(A_i)=D_{i}(A_i)-1$; and for every $A_i < m$, we have $D_i(m)=D_{i+1}(m)$.

For the law of $A_1$, we use equations \eqref{eq:fgpmn} and \eqref{eq:etanueva} to obtain
\begin{align*}
	\mathbb{P}(A_1>n) & = \prod_{i=1}^n \mathbb{P}(\eta^{(i)}=0)  =  \prod_{i=1}^n \frac{1- f_{-i+1,0}(0) }{1- f_{-i,0}(0)} f^\prime_{-i+1} (f_{-i+1,0}(0) )  = \frac{ f^{\prime}_{-n,0}(0) }{1-f_{-n,0}(0)}.
	\end{align*}

\hfill	\qed

The previous techniques can be adapted to establish Theorem \ref{th:ThB}.

\subsection{Proof of Theorem \ref{th:ThB}}
Thanks to \eqref{eq:AdeD} and \eqref{eq:defli}, we establish \eqref{eq:AdeB}. 
Equation \eqref{eq frak a} implies that $D_{i+1}(A_i) = D_i(A_i) -1$ and $D_i(m)=D_{i+1}(m)$ for all $A_i < m \leq \ell$. 
Now, let us analyze the case $1\leq m < A_i$. Define 
	\[
	E_i:=\{ B_j=(d_j(1),\dots, d_j(l_j)): 1\leq j\leq i \},
	\]
	where $\{l_j: 1\leq j\leq i\}$ and $\{d_j(n):1\leq n\leq l_j \, \mbox{ and }\, 1\leq j\leq i \}$ are fixed integer numbers and $l_i=\ell$. Conditional on $E_i$,  $A_i$ only depends on the vector $B_i$. Furthermore, $E_i$ only depends on the genealogical information of the descendants of individual $(-(\ell+1), \mathfrak{a}_i( \ell+1) )$. Then, we need  to consider the tree $\mathcal{T}^{(-(\ell+1), \mathfrak{a}_i( \ell+1) )}$,  the stopping line $L_{\ell+1}$ given by \eqref{def:stoppingline} with $k=\ell+1$, and the set of roots $R_{\ell+1}$ given by \eqref{def:roots} with $k=\ell+1$. Denote by $\mathcal{G}_{\ell+1}$ the $\sigma$-algebra generated by the subtrees with roots in $R_{\ell+1}$. See the left tree on Figure \ref{fig:2} with $k=\ell+1$. 
	
Since $E_i  \in \mathcal{G}_{\ell+1} \vee \mathcal{F}_{L_{\ell+1}} $, according to equation \eqref{eq:cinlar} we have
\begin{align}
	\mathbb{P}(  \mathcal{T}^{(- (A_{i}-1) ,\mathfrak{a}_{i+1}(A_{i}-1)  )} = {\bf t} \mid E_i)
& = \mathbb{P}(  \mathbb{P} ( \mathcal{T}^{( - (A_{i}-1) ,\mathfrak{a}_{i+1}(A_{i}-1)  )} ={\bf t} \mid \mathcal{G}_{\ell+1} \vee \mathcal{F}_{L_{\ell+1}}  )  \mid E_i) \nonumber \\
& = \mathbb{P}(  \mathbb{P} ( \mathcal{T}^{(- (A_{i}-1),\mathfrak{a}_{i+1}(A_{i}-1)  )} = {\bf t} \mid \mathcal{F}_{L_{\ell+1}}  )    \mid E_i).\label{eq:torre law}
\end{align}
Thanks to equation \eqref{eq:medq}, we obtain
\begin{equation*}
 \mathbb{P}(  \mathcal{T}^{(- (A_{i}-1),\mathfrak{a}_{i+1}(A_{i}-1)  )} = {\bf t}  \mid E_i) =  \mathbb{P}(  \mathbb{Q}_{(-A_i+1)}  \left(    {\bf t} \right)   \mid E_i).
\end{equation*}
		Recall that the vector $(D_{i+1}(1), D_{i+1}(2), \dots, D_{i+1}(A_{i} -1) )$ is a functional of the tree  $\mathcal{T}^{(-(A_i-1),\mathfrak{a}_{i+1}(A_{i}-1)  )}$. Therefore, conditional on $E_i$, $(D_{i+1}(1), D_{i+1}(2), \dots, D_{i+1}(A_{i} -1) )$  only depends on the vector $B_i$, because conditional on $E_i$, $A_i$ only depends on $B_i$.

Continuing with the same arguments as in the proof of Theorem \ref{th:ThB}, we get \eqref{eq:same}. From this equation, it follows that for every $1\leq m <A_i$, the variable $D_{i+1}(m)$ conditioned on $E_i$  has the same  distribution as $\eta^{(-m)}$  given in \eqref{eq:etadistribucion}. Furthermore, $D_{i+1}(A_i)=D_{i}(A_i)-1$; and for every $A_i < m \leq \ell$, we have $D_i(m)=D_{i+1}(m)$. Therefore, the vector $(D_{i+1}(1), \dots, D_{i+1}(\ell) )$ given $E_i$ only depends on the values of $B_i$.
		
If $D_{i+1}(m)\neq 0$ for  some $1\leq m\leq \ell$, by \eqref{eq:AdeB} we have $A_{i+1}\leq \ell$. Recalling equation \eqref{eq:defli}, we conclude that $l(B_{i+1})=l(B_i)\vee A_{i+1}=\ell $ and $B_{i+1}=(D_{i+1}(m), 1\leq m\leq \ell)$. Therefore, we have determined the distribution of $B_{i+1}$ conditioned on $E_{i}$. 
		Otherwise, if $D_{i+1}(m)= 0$ for every $1\leq m\leq \ell$. By \eqref{eq:defli} and \eqref{eq:AdeB}, we deduce that  $A_{i+1}> \ell$ and $l(B_{i+1})>\ell$. Now, by the definition of $B_{i+1}$, we need to determine the law of  $l(B_{i+1})$  and  $\{D_{i+1}(m), \ell<m\leq l(B_{i+1})\}$ conditioned on $E_i$. To achieve this, we analyze the tree $\mathcal{T}^{(-( \ell+j) , \mathfrak{a}_i( \ell+j) )}$ for a fixed $j\geq1$.
		 Consider the stopping line $L_{\ell+j}$ defined by \eqref{def:stoppingline}. As previously, $\mathcal{F}_{L_{\ell+j}}=\sigma(\Omega_{(-n,\mathfrak{a}_{i}(n))}: A_i \leq  n\leq \ell +j )$. 
	Observe that $E_i  \in \mathcal{G}_{\ell+1} \vee \mathcal{F}_{L_{\ell+j}} $, where $\mathcal{G}_{\ell+1} $ is the $\sigma$-algebra generated by the subtrees with roots in $R_{\ell+1}$.    
According to Proposition \ref{pr:lineaparo}, given $\mathcal{F}_{L_{\ell+j}}$, the subtrees $\mathcal{T}^{(x)}$ with roots $x\in R_{\ell+1}$ and the subtrees $\mathcal{T}^{(y)}$ with roots $y\in L_{\ell+j}\setminus L_{\ell+1}$ are independent. See the right tree on Figure \ref{fig:2}.
Therefore,  employing the same argument as in \eqref{eq:cinlar}, for every  $y\in L_{\ell+j}\setminus L_{\ell+1}$, we have
\begin{equation*}
\mathbb{P}\left(  \mathcal{T}^{(y)} = {\bf t}\ \Big|\  E_i\right)
= \mathbb{P}(  \mathbb{P} ( \mathcal{T}^{( y )} ={\bf t} \mid \mathcal{G}_{\ell+1} \vee \mathcal{F}_{L_{\ell+j}} )  \mid E_i) 
 = \mathbb{P}\left(  \mathbb{P} ( \mathcal{T}^{(y)} = {\bf t} \mid   \mathcal{F}_{L_{\ell+j}}  )    \ \Big| \ E_i\right).
\end{equation*}
Now, if $y\in L_{\ell+j}\setminus L_{\ell+1}$, there exists an $1\leq r<j$ such that $|y|=r$, where $|y|$ denotes its length in  $\mathcal{T}^{(-( \ell+j) , \mathfrak{a}_i( \ell+j) )}$.  We utilize \eqref{eq: propapl} to obtain 
\begin{equation*}
\mathbb{P} ( \mathcal{T}^{(y )} = {\bf t} \mid   \mathcal{F}_{L_{\ell+j}}  ) =\mathbb{Q}_{(-(\ell+j)+r)}  \left(    {\bf t} \right).
\end{equation*}
Due to the fact that the right-hand side of the previous equation, conditioned on $E_i$, only depends on $\ell$, we obtain 
\begin{equation}\label{ref:moved}
\mathbb{P}\left(  \mathcal{T}^{(y)} = {\bf t}\ \Big|\  E_i\right) = \mathbb{Q}_{(-(\ell+j)+r)}  \left(    {\bf t} \right).
\end{equation}
In other words, given $E_i$,  every subtree with root in $y\in L_{\ell+j}\setminus L_{\ell+1}$ remains a GWVE tree with environment $\mathcal{S}_{(\ell+j-|y|)}$. 
Recall that $D_{i+1}(m)= | \mathcal{D}_{i+1}(m) | -1$, where 
\[
\mathcal{D}_{i+1}(m)
		=  \{ \text{daughters of } (-m,\mathfrak{a}_{i+1}(m)) \text{ with descendants in } \{ (0,k): k\geq i+1  \} .
\]
Moreover, $\mathfrak{a}_{i}(m)=\mathfrak{a}_{i+1}(m)$ for $m\geq \ell$, and  $\ell$ is the coalescent time of individuals  $\{(0,k): k\leq i+1 \}$. Consequently, for every $m\geq \ell+1$, given $E_i$,  $\{(0,k): k\leq i+1 \}$  are descendant of  $(-m, \mathfrak{a}_{i}(m) )$ and 
\[
\mathcal{D}_{i+1}(m)
		= \{ \text{daughters of } (-m,\mathfrak{a}_{i}(m)) \text{ with descendants in } \{ (0,k): k\geq 1  \}  \}.
\]
 Now, let us examine the tree $\mathcal{T}^{(-( \ell+j) , \mathfrak{a}_i( \ell+j) )}$. Note that the daughters of individual $(-m,\mathfrak{a}_{i}(m))$ have length $r=\ell+j-m+1$. Furthermore, $(-(m-1),\mathfrak{a}_{i}(m-1)) $  is a daughter of $(-m,\mathfrak{a}_{i}(m)) $ whose descendants are alive in the present generation. Thus, in the definition of $D_{i+1}(m)$, the subtraction one can be  attributed to individual 
$(-(m-1),\mathfrak{a}_{i}(m-1)) $, and for every $\ell+1\leq m<\ell+j$, the value  $D_{i+1}(m)$  can regarded as a functional of the subtrees $\mathcal{T}^{(y)}$ with $y\in L_{\ell+j}\setminus L_{\ell+1}$ such that $|y|=\ell+j-m+1$. By \eqref{ref:moved}, we know that these subtrees, conditioned on $E_i$, follow the law given by $\mathbb{Q}_{(-m+1)}$. Consequently, the distribution of $D_{i+1}(\ell+m)$ conditioned on $E_i$ is determined by \eqref{eq:etadistribucion}, for every $\ell+1\leq m<\ell+ j$. Since $j$ was fixed but arbitrary, we have shown that for every $\ell+1\leq m$, the distribution of $D_{i+1}(\ell+m)$ conditioned on $E_i$ only depends  on $\ell$ and is given by \eqref{eq:etadistribucion}.
	
Recall that by the definition of $B_{i+1}$, its length  $l(B_{i+1})$ is the first $n$ such that $D_{i+1}(n)\neq 0$. Consequently,  $l(B_{i+1})$ conditioned on $E_i$ follows the same distribution as \eqref{eq:LMarkov}, and we have established the Markov property with its transitions.

\hfill	\qed

\subsection{Proof of Proposition \ref{pr:Prop Linear}}
	Observe that for every $n\geq 1$, the variable $\eta^{(-n)}\overset{\mathcal{L}}{=}\eta_n$ with environment $\mathcal{E}_n=( e_{-n+1}, e_{-n+2}, \dots , e_{0})$ as shown in Example \ref{Ex linear1}. In particular, $\eta^{(-n)}$ follows a geometric distribution (modeling the number of failures until the first success) with success probability $\lambda_1:=p_0$ and for $n\geq 2$,
	\begin{equation}\label{eq:failure prob}
	\lambda_n:=\left( 1 + \underset{k=-n+2}{\overset{0}{\sum}}
	s_{k}   \right)\left( 1 + \underset{k=-n+1}{\overset{0}{\sum}}
	s_k   \right)^{-1}.
	\end{equation}
	where $s_{0}=\tfrac{1-p_0}{p_0}$ and $s_{k}=\tfrac { (1-p_k) }{p_{k}}\tfrac{r_{k+1}  \cdots r_0 }{  p_{k+1}\cdots p_0}\,$ for $k<0$.
	By induction on $i\geq 1$, we are going to prove the following statement:
	\begin{description}
		\item[(${H_i}$)]\label{H_i} The random variables $(D_i(n),n\geq 1)$ are independent, geometrically distributed with success probability $\lambda_n$ given in \eqref{eq:failure prob}. Additionally, they are independent of $(A_0,\dots, A_{i-1})$.
	\end{description}
	The claim holds once we prove that $(H_i)$ is true for every $i\geq 1$. Indeed, suppose that $(H_i)$ is true for every $i\geq 1$. Then, by equation \eqref{eq:AdeD}, we can see that $A_i$ is independent of $(A_0,\dots, A_{i-1})$ and it is distributed as $A_1$. In particular,
	\begin{equation*}
	\mathbb{P}(A_i>n)= \prod_{\ell=1}^n \mathbb{P}(\eta^{(\ell)}=0)= \left( 1 + \underset{k=-n+1}{\overset{0}{\sum}}
	s_k   \right)^{-1}, \qquad n\geq 1.
	\end{equation*}
	
Observe that $(H_1)$ holds by Proposition \ref{pr:ThD}. Now, assuming $(H_i)$, we will prove that $(H_{i+1})$ is also true by conditioning on the value of $A_i$. Suppose that $A_i=h$. We apply the transition probabilities from Proposition \ref{pr:ThD}. Note that  $D_i(n)=D_{i+1}(n)$ for all $n>h$. Then, by hypothesis $(H_i)$, the variables $(D_{i+1}(n), n>h)$  are independent, geometrically distributed with parameters $\lambda_n$ and also independent of $(A_0,\dots, A_{i-1})$. For $n=h$, by hypothesis, $D_{i+1}(h)=D_i(h)-1$ is independent of $(D_i(n),n> h)$ and independent 
	of $(A_0,\dots, A_{i-1})$. Additionally, by \eqref{eq:AdeD}, $D_i(h)$ is a geometric variable with parameter $\lambda_h$  conditioned to be strictly positive. Then,
	\[
	\mathbb{P}(D_{i+1}(h)=k)=\frac{\mathbb{P}(\eta^{(h)}=k+1)}{\mathbb{P}(\eta^{(h)}>0)}=\lambda_h(1-\lambda_h)^k,\qquad k\geq 0,
	\]
	which is a geometric random variable with parameter $\lambda_h$.
	Finally, $(D_{i+1}(n), n<h)$ are \textit{new} independent geometric random variables with parameters $\lambda_n$. Therefore, they are  independent of $(D_{i+1}(n), n\geq h)$ and $(A_0,\dots, A_{i-1})$. 
	In other words, conditionally on  $\{A_i=h\}$, the variables $(D_{i+1}(n), n\geq 1)$ are independent with geometric distributions of parameters $\lambda_n$ and independent of $(A_0,\dots, A_{i-1})$.  Integrating over $h$ yields  the result $(H_{i+1})$.
\hfill \qed

\section*{Acknowledgments}
The authors express their gratitude to Amaury Lambert and Lea Popovic for their valuable discussions and insightful comments.  S. P.  was supported by UNAM-DGAPA-PAPIIT grant no. IA103220.

\end{document}